\documentclass[a4paper,12pt]{article}
\usepackage{amsmath,amssymb}
\usepackage{enumerate}
\usepackage{epic}
\usepackage{bbm}
\usepackage{theorem}

\usepackage{bbm}
\usepackage{mathptmx}

\newcommand{\NN}{\mathbb{N}}
\newcommand{\ZZ}{\mathbb{Z}}

\newcommand{\cL}{\mathcal{L}}
\newcommand{\cB}{\mathcal{B}}

\newcommand{\<}{\ensuremath{\langle}}
\renewcommand{\>}{\ensuremath{\rangle}}

\newcommand{\ini}{\operatorname{in}}
\newcommand{\bin}{\operatorname{bin}}

\newcommand{\sat}{\operatorname{sat}}

\newcommand{\Span}{\operatorname{span}}

\renewcommand{\>}{\rangle}

\theorembodyfont{\rmfamily}
\newtheorem{Theorem}{Theorem}[section]
\newtheorem{Algorithm}[Theorem]{Algorithm}

\newtheorem{Proposition}[Theorem]{Proposition}

\newtheorem{Example}[Theorem]{Example}
\newtheorem{Remark}[Theorem]{Remark}

\begin{document}
\title{Homogeneous Buchberger algorithms and Sullivant's
computational commutative algebra challenge}
\author{
Niels Lauritzen\\
Institut for Matematiske Fag\\
Aarhus Universitet\\
DK-8000 \AA rhus\\
Denmark\\ 
niels@imf.au.dk
}

\maketitle

\begin{abstract}
We give a variant
of the homogeneous Buchberger algorithm for positively graded lattice ideals.
Using this algorithm we solve 
the Sullivant computational commutative 
algebra challenge\footnote{
{\tt http://math.berkeley.edu/\~{}seths/ccachallenge.html}}.
\end{abstract}


\section{Introduction}
Suppose that $I$ is a homogeneous ideal in
a polynomial ring 
$R = k[x_0, \dots, x_n]$ over a field $k$. The usual
homogeneous Buchberger algorithm builds a Gr\"obner basis
for $I$ by successively constructing
truncated Gr\"obner bases of increasing degrees.
Suppose that $I$ is {\it saturated\/} i.e. $I = \bar{I} = \{g\in R\mid
(x_0\cdots x_n)^m g \in I, m\gg 0\}$.
If we encounter a polynomial $f$
divisible by a variable in degree $d$ of the
homogeneous Buchberger algorithm, then we may conclude that
$f$ reduces to zero modulo the already constructed
truncated Gr\"obner basis in degree $<d$ for $I$. 
This simple observation also allows for detection
of non-saturated ideals in some cases.

Sullivant's challenge is about deciding if a specified set $B$
of $145,512$ binomials generate the kernel $P$ of the (toric)
ring homomorphism
$\varphi: k[x_{ijk}] \rightarrow k[u_{ij}, v_{ik}, w_{jk}]$ given by 
$$
\varphi(x_{ijk}) = u_{ij} v_{ik} w_{jk},
$$
where $1\leq i,j,k \leq 4$. 
We give a version of the homogeneous Buchberger algorithm
with a Gebauer-M\"oller criterion specifically tailored to
positively graded lattice ideals. 
Using an implementation of this algorithm in the software
package {\tt GLATWALK}\footnote{\tt http://home.imf.au.dk/niels/GLATWALK} 
we deduce that {\it the ideal $J$ generated by $B$ is strictly contained in 
$P$\/} by showing that $J$ cannot be saturated. In fact, we exhibit a specific
binomial $b$ of degree $14$ in $\bar{J}\setminus J$.

I am grateful to B.~Sturmfels for stimulating my interest in 
Sullivant's computational commutative algebra
challenge. R.~Hemmecke has made me aware that he and P.~Malkin
already computed the full Gr\"obner basis of $J$ using
new algorithms in a new version of {\tt 4ti2} thereby answering Sullivant's
challenge. In fact they prove that
the ``missing'' binomials in Sullivant's challenge have degree $14$ and 
form an orbit under the action of a certain symmetry group.
I am grateful to Hemmecke for verifying that $b$ lies in this orbit.

\section{Preliminaries}

We let $R = k[x_1, \dots, x_n]$ denote the ring of polynomials
over a field $k$. We assign degrees to the variables by 
$\deg(x_1) = a_1, \dots, \deg(x_n) = a_n$, where $a_1, \dots,
a_n$ are positive integers. A monomial $x^v\in R$ has degree $\deg(x^v) = 
v_1 a_1 + \cdots + v_n a_n$, where $v = (v_1, \dots, v_n)$.
This gives the (positive) grading
$$
R = \oplus_{s\geq 0} R_s,
$$
where $R_s = \Span_k \{x^v \mid \deg(x^v) = s\}$. For a monomial
order $\prec$ on $R$ and a subset $S\subset R$ we let
$\ini_\prec(S) = \{\ini_\prec(f) \mid f\in S\setminus\{0\}\}$. A
Gr\"obner basis for an ideal $I\subset R$ is a finite subset
$G\subset I$, such that $\<\ini_\prec(G)\> = \ini_\prec(I)$.

\subsection{Truncated Gr\"obner bases}

For a homogeneous ideal $I$ in $R$ and $d\in \NN$ we let
$$
I_{<d} = \bigoplus_{s< d} I_s.
$$
A $d$-truncated Gr\"obner basis for $I$ is a finite subset
$G_{<d}\subset I_{<d}$, such that $\<\ini_\prec(G_{<d})\>_{<d} 
= \ini(I)_{<d}$ i.e.
we require only match of initial ideals up to degree $d$. Using
the division algorithm it is easy to show that $f\in I_{<d}$ reduces
to zero modulo the polynomials in a $d$-truncated Gr\"obner basis for
$I$.

\section{The homogeneous Buchberger algorithm with sat-reduction}

We call an ideal $I$ {\it saturated\/} if $I= \bar{I} =
\{g\in R\mid (x_0\cdots x_n)^m g \in I, m\gg 0\}$. This means
that $m f\in I$ implies $f\in I$, where $m$ is a monomial and
$f$ a polynomial in $R$.
Let $\prec$ be a term order on $R$. 
For a polynomial $f\in R$ we let
$\sat(f)$ denote $f$ divided by the greatest common divisor of the
monomials in $f$. We say
that $f$ {\it sat-reduces\/} to $h$ modulo $g$ if either $h = \sat(f)$ and 
$\deg(h) < \deg(f)$ or $f$ reduces to $h$ modulo $g$ in the
usual sense i.e. 
$\ini_\prec(g)$ divides
a term $t$ in $f$ and
$$
h = f - (t/\ini_\prec(g)) g.
$$
Notice that
if $f$ sat-reduces to $h$ modulo $g$ and $f, g$ belongs to 
a saturated ideal $I$, then
$h\in I$.
A remainder in the division algorithm of $f$ modulo a set
of polynomials $G$ using
sat-reduction in each step is denoted $f^{G(\sat)}$. 

The $S$-polynomial of two homogeneous polynomials
is homogeneous of degree no less than the degrees of the polynomials.
The (usual) reduction of a homogeneous polynomial of degree $d$ modulo
a set of homogeneous polynomials gives a homogeneous polynomial
of degree $d$.  These observations give the 
homogeneous Buchberger
algorithm as explained in (\cite{CKR}, Theorem 11). We tailor the 
homogeneous Buchberger algorithm to the special case where input consists of a
set $B = \{f_1, \dots, f_r\}\subset R$ of
homogeneous polynomials generating a saturated ideal.
This has the consequence that reduction of a 
homogenous polynomial $f$ of degree $d$ 
divisible by a variable $x_i$ 
is not necessary, since $f/x_i\in I_{<d}$ reduces to zero using the
already computed $d$-truncated Gr\"obner basis $G_{<d}$.

\begin{Algorithm}[Homogeneous Buchberger algorithm for saturated ideals]
\label{AlgorithmBAG}
\

\

\noindent
{\bf INPUT}:  Term order $\prec$. Homogeneous polynomials 
$B = \{f_1, \dots, f_r\} \subset R$ generating a 
saturated homogeneous ideal
$I$.

\

\noindent
{\bf OUTPUT}: Homogeneous polynomials $G = \{g_1, \dots, g_s\}$ such
that $\{g_1, \dots, g_s\}$ is a minimal Gr\"obner
basis over $\prec$ for the ideal generated by $B$.

\begin{enumerate}[(i)]
\item $Spairs := \emptyset$; $G := \emptyset$;
\item while ($B\neq \emptyset$ or $Spairs \neq \emptyset$) do
\begin{enumerate}
\item Extract\footnote{This means that $f$ is deleted from the
    relevant list after it is extracted}  
a polynomial $f$ of minimal degree in $B\cup Spairs$. 
\item Compute $g := f^{G(\sat)}$, continue if the degree drops in a
sat-reduction step in the division algorithm;\label{Stepsat}
\item if ($g = 0$) continue;
\item \label{Steplast} $G := G \cup \{g\}$;
\item \label{Sstep} Append $S$-polynomials $S(g, h)$ to $Spairs$ for every $h\in G\setminus\{g\}$.

\end{enumerate}
\end{enumerate}

\end{Algorithm}

\begin{Remark}

\

\begin{enumerate}[(i)]
\item
After step (\ref{Steplast}) in Algorithm \ref{AlgorithmBAG}, 
the polynomials of degree $<d$ in $G$ form a minimal $d$-truncated 
Gr\"obner basis of $I$, where $d$ is the minimal degree of 
the polynomials in $B\cup Spairs$.
\end{enumerate} 
\end{Remark}

An easy modification to algorithm (\ref{AlgorithmBAG}) may detect if
$I$ is not saturated. If the sat-reduction $f^{G(\sat)}$ 
of $f$ is non-zero and has lower degree than $f$, then we may deduce
the existence of a monomial $x^v$ and a polynomial $g$ such that
$x^v g\in I$, but $g\not\in I$.

\begin{Algorithm}[Homogeneous Buchberger algorithm with saturation check]\label{AlgorithmBAGcheck}
\

\

\noindent
{\bf INPUT}:  Homogeneous polynomials 
$B = \{f_1, \dots, f_r\}$ and a term order $\prec$.

\

\noindent
{\bf OUTPUT}: Homogeneous polynomials $G = \{g_1, \dots, g_s\}$ such
that $\{g_1, \dots, g_s\}$ is a minimal Gr\"obner
basis over $\prec$ for the ideal $I$ generated by $B$ or proof that
$I$ is not saturated.

\begin{enumerate}[(i)]
\item $Spairs := \emptyset$; $G := \emptyset$;
\item while ($B\neq \emptyset$ or $Spairs \neq \emptyset$) do
\begin{enumerate}
\item Extract a polynomial $f$ of minimal degree $d$ in $B\cup Spairs$. 
\item $g := f^{G(\sat)}$;\label{Stepr}
\item if ($g = 0$) continue;
\item if($\deg(g) < d$)
\begin{enumerate}[(i)]
\item OUTPUT $f$ as proof that $I$ is not saturated and HALT.
\end{enumerate}
\item $G := G \cup \{g\}$;
\item \label{Sstepcheck} Append $S$-polynomials $S(g, h)$ to $Spairs$ for every $h\in G\setminus\{g\}$.
\end{enumerate}
\end{enumerate}

\end{Algorithm}

\begin{Example}
We give a simple example illustrating algorithm (\ref{AlgorithmBAGcheck}).
\begin{enumerate}[(i)]
\item
Consider the input $B = \{x z - y^2, x^4 - y^3\}$ along with
the reverse lexicographic term order $x\prec y \prec z$. 
\item
The ideal $I$ generated
by $B$ is homogeneous in the grading $\deg(x) = 3, \deg(y) = 4,
\deg(z) = 5$ and $\deg(x z - y^2) = 8 < \deg(x^4 - y^3) = 12$. 
\item
After the first
loop we have $B = \{y^3 - x^4\}, G = \{y^2 - x z\}$ and $Spairs = \emptyset$,
where $G$ is a $12$-truncated Gr\"obner basis of $I$.
\item
In the second loop we sat-reduce $y^3 - x^4$ modulo $y^2 - x z$ and
get $y z - x^3$. As $\deg(y z - x^3) = 9 < \deg(y^3 - x^4) = 12$, we
conclude that $I$ is not saturated.
\end{enumerate}

\begin{enumerate}[(i)]
\item
Now suppose that $B = \{y^2 - x z, y z - x^3\}$ in
the same grading.
\item
After the second loop we have
\begin{align*}
B &= \emptyset\\
G &= \{y^2 - x z, y z - x^3\}\\
Spairs &= \{y x^3 - z^2 x\},
\end{align*}
where $G$ is a $13$-truncated Gr\"obner basis
of $I$.
\item
Now $y x^3 - z^2 x$ sat-reduces to $z^2 - y x^2$ modulo $G$. We conclude
that $I$ is not saturated.
\end{enumerate}
\begin{enumerate}[(i)]
\item
Now proceed with $B=\{y^2 - x z, y z - x^3, z^2 - y x^2\}$.
\item
After a few loops we have
\begin{align*}
B &= \emptyset\\
G &= \{y^2 - x z, y z - x^3, z^2 - y x^2\}\\
Spairs &= \{y^2 x^2 - z x^3\},
\end{align*}
where $G$ is a $14$-truncated Gr\"obner basis
of $I$. Since $y^2 x^2 - z x^3$ sat-reduces to zero, $G$ is the
reduced Gr\"obner basis of $I$.
\end{enumerate}
\end{Example}

The number of $S$-pairs considered for reduction can be reduced drastically
by using a version of the Gebauer-M\"oller criterion in algorithms 
(\ref{AlgorithmBAG}) and (\ref{AlgorithmBAGcheck}). 
The framework for properly explaining 
the Gebauer-M\"oller criterion is in the 
context of Gr\"obner bases for modules (cf. \cite{CKR}, \S4).

\subsection{The Gebauer-M\"oller criterion}

Let $e_1, \dots, e_m$ denote the canonical basis of the finitely
generated free module $F = R^m$. A monomial in $F$ is an
element $x^v e_i$, where $x^v$ is a monomial in $R$. 
Every element in $F$ is a 
$k$-linear combination of monomials.
By definition a monomial
$x^\alpha e_i$ divides a monomial $x^\beta e_j$ if and only if
$i = j$ and $x^\alpha$ divides $x^\beta$ in $R$. We write this as
$x^\alpha e_i \mid x^\beta e_j$. A monomial order on $F$ is
a total order $\prec$ on monomials in $F$ satisfying
$$
x^\alpha e_i \prec x^\beta e_j \implies x^{\alpha + \gamma} e_i \prec 
x^{\beta + \gamma} e_j
$$
for every $i, j = 1, \dots, m$ and $\alpha, \beta, \gamma\in \NN^n$. 
We let $\ini_\prec(f)$ denote
the largest monomial in $f$. Now the Gr\"obner basics for ideals in $R$
can be generalized to submodules of $F$ almost verbatim. For a
subset $B\subset F$ we let $\ini_\prec(B)$ denote the submodule generated
by $\ini_\prec(f)$, where $f\in B$.
A Gr\"obner
basis of a submodule $M \subset F$ is a set of elements
$G = \{m_1, \dots, m_t\} \subset M$ satisfying $\ini_\prec(M) =
\ini_\prec(G)$.
It is called minimal if $\ini_\prec(m_i) \nmid \ini_\prec(m_j)$ for $i\neq j$.
We will use Gr\"obner bases for submodules in reasoning about syzygies
of monomial ideals. Consider a monomial ideal 
$$
M = \<x^{v_1}, \dots, x^{v_m}\>\subset R.
$$ 
The syzygies of $M$ are the relations in $M$ i.e. the
kernel $K$ of the natural surjection $R^m \rightarrow M$.
Now consider the $\ZZ^n$-grading $\deg(x_i) = e_i$ on $R$. 
Then $K$ is a homogeneous submodule of $F$ in the $\ZZ^n$-grading given by
$\deg(e_i) = v_i$. A natural set of homogeneous
generators are
$$
S_{ij} = x^{v_i \vee v_j - v_j} e_j - x^{v_i\vee v_j - v_i} e_i\in K
$$
for $1\leq i < j \leq m$ (see \cite{CLO}, Proposition 2.8). 
Define a monomial order $\prec$ (The Schreyer order) 
on $F$ by $x^\alpha e_i \prec x^\beta e_j$ if and only if
$$
\alpha + v_i < \beta + v_j\mbox{\ or\ }\alpha + v_i = \beta + v_j
\mbox{\ and\ }i < j,
$$ 
where $<$ is any term order on $R$. Then we have
the following
\begin{Proposition}\label{PropositionGM}
The homogeneous generating set $\{S_{ij} \mid 1\leq i < j \leq m\}$ 
is a Gr\"obner basis for $K$ over the Schreyer order $\prec$.
\end{Proposition}

The Gr\"obner basis in Proposition \ref{PropositionGM} is rarely
minimal. In view of Theorem 2.9.9 in \cite{CLO}, it suffices to reduce
the $S$-pairs corresponding to a minimal Gr\"obner basis of
the syzygies (in Buchberger's algorithm). This procedure
is in fact one of the Gebauer-M\"oller criteria for cutting down
on the number of $S$-pairs.
The point is that this minimization is
easy and quite fast to perform in step (\ref{Sstep}) of 
Algorithm \ref{AlgorithmBAG}. Suppose that we must update
$Spairs$ with a non-zero polynomial $g = g_m$, where
$G = \{g_1, \dots, g_{m-1}\}$ in step (\ref{Sstep}). We put
$x^{v_i} = \ini_\prec(g_i)$ for $i = 1, \dots, m$.
Consider the syzygies $S_{1m}, \dots, S_{m-1, m}$. In the Schreyer 
order we have $\ini_\prec(S_{1m}) = x^{v_1\vee v_m - v_m} e_m,
\dots, \ini_\prec(S_{m-1, m}) = x^{v_{m-1}\vee v_m - v_m} e_m.$
Thus the minimization can be done successively in step (\ref{Sstep}) by
throwing out superfluous monomials among 
\begin{align*}
&x^{v_1\vee v_m - v_m}\\
&\vdots\\
&x^{v_{m-1}\vee v_m - v_m}.
\end{align*}
This can be implemented as below ($u\leq v$ means that 
$v-u\in \NN^n$ for $u, v\in \NN^n$), where (\ref{stepdisj})
represents the usual criterion, where leading terms are relatively
prime (cf.~\cite{CLO}, Proposition 2.9.4).

\begin{Algorithm}\label{AlgorithmGM}

\

\

\noindent
{\bf updateSpairs:}
\begin{enumerate}[(i)]
\item
$MinSyz := \emptyset$;
\item
for each $v_i$ in $\{v_1, \dots, v_{m-1}\}$ do
\begin{enumerate}
\item\label{stepdisj}
if ($v_m \wedge v_i = 0$) continue;
\item
$a = v_i \vee v_m - v_m$;
\item
if ($w \leq a$ for some $(w, p)\in MinSyz$) continue;
\item
Delete $(w, p)\in MinSyz$ if $a\leq w$;
\item
$MinSyz := MinSyz \cup \{(a, S(g_i, g_m))\}$;
\end{enumerate}
\item
for each $(a, p)\in MinSyz$ do 
\begin{enumerate}
\item
$Spairs := Spairs\cup \{p\}$;
\end{enumerate}
\end{enumerate}

\end{Algorithm}

\section{Lattice ideals}

Recall the decomposition of an integer vector 
$v\in \ZZ^n$ into $v = v^+ - v^-$,
where $v^+, v^-\in \NN^n$ are vectors with disjoint support.
For $u, v\in \NN^n$ we let $u \leq v$ denote the partial
order given by $v - u\in \NN^n$.
For a subset
$\cB \subset \ZZ^n$ we associate the ideal
$$
I_\cB = \<x^{v^+} - x^{v^-} \mid v\in \cB\> \subset R.
$$
In the case where $\cB = \cL$ is a lattice we call
$I_\cL$ the {\it lattice ideal\/} associated to $\cL$. 
If $u - v\in \cL$ for $u, v\in \NN^n$, then
\begin{equation}\label{satref}
x^u - x^v = x^{u - (u - v)^+} (x^{(u-v)^+} - x^{(u-v)^-})\in I_\cL.
\end{equation}
The binomials $B_\cL = \{x^u - x^v\mid u-v\in \cL\}\subset I_\cL$ 
are stable under the fundamental operations in
Buchberger's algorithm: forming $S$-polynomials and reducing modulo
a subset of $B_\cL$. This means that starting with a generating set
for $I_\cL$ in $B_\cL$ we end up with a Gr\"obner basis consisting
of binomials in $B_\cL$. Reducing a monomial $x^w$ by an element
of $B_\cL$ amounts to replacing $x^w$ by $x^{w - v}$, where $v\in \cL$.
Therefore if a binomial
$x^u - x^v\in I_\cL$, then $u - v \in \cL$. This proves that
$I_\cL$ is saturated and algorithm 
(\ref{AlgorithmBAG}) applies. The simple data structures in the specialization
of algorithm (\ref{AlgorithmBAG}) to lattice ideals are very appealing.
If $f = x^u - x^v$, then
$$
\sat(f) = x^{(u-v)^+} - x^{(u-v)^-}.
$$
by (\ref{satref}).
With this in mind we define 
$$
\bin(w) = x^{w^+} - x^{w^-}
$$
for $w\in \ZZ^n$. Using this notation we have $\sat(\bin(u), \bin(v)) = 
\bin(u-v)$. Similarly if 
$v^+ \leq u^+$ we may reduce $\bin(u)$ by
$\bin(v)$. This results in a binomial $f$ with $\sat(f) = \bin(u-v)$. 
Notice that replacing $u$ by $u - v$ if $v^+ \leq u^+$ corresponds
to sat-reduction of $\bin(u)$ by $\bin(v)$.
We have silently assumed that
the initial term of $\bin(w)$ is $x^{w^+}$ for the term order
in question. We will keep this convention throughout.

Usually a generating set $\cB$ for $\cL$ as an abelian group is given. 
Computing the lattice ideal $I_\cL \supset I_\cB$ can be done using that
$$
I_\cL = \bar{I}_{\cB}.
$$
If $\cB$ contains a positive vector, then $I_\cB = I_\cL$ (\cite{St}, 
Lemma 12.4). If $\cL \cap \NN^n = \{0\}$, 
$I_\cL$ may be computed from $I_\cB$ using 
Gr\"obner basis computations for different reverse 
lexicographic term orders (\cite{St}, Lemma 12.1).

With these conventions it is quite easy to convert algorithm
(\ref{AlgorithmBAG}) into a specialized algorithm for
lattice ideals representing binomials via integer
vectors with additional structure (like the degree of $\bin(v)$
and certain other (optimizing) features). We give the
straightforward translation of algorithm (\ref{AlgorithmBAG}) 
into the lattice case.

\begin{Algorithm}[Homogeneous Buchberger algorithm for lattice ideals]\label{AlgorithmBA}
\

\

\noindent
{\bf INPUT}: Term order $\prec$. Integer vectors 
$B = \{v_1, \dots, v_r\}$ with respect
to $\prec$ such that $\<\bin(v_1), \dots, \bin(v_r)\>$ is a 
positively graded lattice ideal $I_\cL$.

\

\noindent
{\bf OUTPUT}: Integer vectors $G = \{w_1, \dots, w_s\}$ such
that $\<\bin(w_1), \dots, \bin(w_s)\>$ is a minimal Gr\"obner
basis over $\prec$ for $I_\cL$.

\begin{enumerate}[(i)]
\item $Spairs := \emptyset$; $G := \emptyset$;
\item while ($B\neq \emptyset$ or $Spairs \neq \emptyset$) do
\begin{enumerate}
\item Extract a binomial $\bin(v)$ of minimal degree in $B\cup Spairs$. 
\item Compute the reduction $\bin(w) := \bin(v)^{G(\sat)}$, continue if the degree drops
in a sat-reduction step in the division algorithm.
\item if ($\bin(w) = 0$) continue;
\item $G := G \cup \{\bin(w)\}$;
\item {\bf updateSpairs} 
\end{enumerate}
\end{enumerate}

\

\noindent
{\bf updateSpairs:}
\begin{enumerate}[(i)]
\item
$MinSyz := \emptyset$;
\item
for each $\bin(v)$ in $G\setminus\{\bin(w)\}$ do
\begin{enumerate}
\item
if ($w^+ \wedge v^+ = 0$) continue;
\item
$a = v^+ \vee w^+ - w^+$;
\item
if ($u \leq a$ for some $(u, p)\in MinSyz$) continue;
\item
Delete $(u, p)\in MinSyz$ if $a\leq u$;
\item
$MinSyz := MinSyz \cup \{(a, \bin(u - v)\}$;
\end{enumerate}
\item
for each $(a, \bin(u))\in MinSyz$ do 
\begin{enumerate}
\item
$Spairs := Spairs\cup \{\bin(u)\}$;
\end{enumerate}
\end{enumerate}

\

\end{Algorithm}

Similarly algorithm \ref{AlgorithmBAGcheck} translates into

\begin{Algorithm}\label{AlgorithmBAcheck}

\

\

\noindent
{\bf INPUT}:  
Term order $\prec$. Normalized integer vectors 
$B = \{v_1, \dots, v_r\}$ with respect
to $\prec$, such that $\<\bin(v_1), \dots, \bin(v_r)\>$ generates 
the ideal $I$.

\

\noindent
{\bf OUTPUT}: 
Integer vectors $G = \{w_1, \dots, w_s\}$ such
that $\<\bin(w_1), \dots, \bin(w_s)\>$ is a minimal Gr\"obner
basis over $\prec$ for $I$ or proof that $I$ is not a lattice ideal.

\begin{enumerate}[(i)]
\item $Spairs := \emptyset$; $G := \emptyset$;
\item while ($B\neq \emptyset$ or $Spairs \neq \emptyset$) do
\begin{enumerate}
\item Extract a binomial $\bin(v)$ of minimal degree $d$ in $B\cup Spairs$. 
\item $\bin(w) := \bin(v)^{G(\sat)}$;
\item if ($\bin(w) = 0$) continue;
\item if($\deg(\bin(w)) < d$)\label{Stephalt}
\begin{enumerate}[(i)]
\item OUTPUT $\bin(w)$ as proof that $I$ is not a lattice ideal and HALT.
\end{enumerate}
\item $G := G \cup \{\bin(w)\}$;
\item {\bf updateSpairs} 
\end{enumerate}
\end{enumerate}

\end{Algorithm}

\section{The Sullivant challenge}

Sullivant's challenge\footnote{\tt
  http://math.berkeley.edu/\~{}seths/ccachallenge.html} is about
deciding if the ideal $J$ 
generated by a given set $B$ of
$145,512$ binomials generate
the kernel $P$ of the toric ring homomorphism 
$$
k[x_{ijk} ] \rightarrow k[u_{ij}, v_{ik}, w_{jk}]
$$
given by $x_{ijk}\mapsto u_{ij} v_{ik} w_{jk}$, where $1\leq i, j, k\leq 4$.
The $145,512$ binomials are constructed by acting with a symmetry group
on carefully selected binomials\footnote{See 
{\tt http://math.berkeley.edu/\~{}seths/ccachallenge.ps} for details}.
In this setting we need to compute in the polynomial ring $k[x_{ijk}]$
in $64$ variables! The ideal $J$ is homogeneous in
the natural grading $\deg(x_{111}) = \cdots = \deg(x_{444}) = 1$.
The strategy is applying algorithm (\ref{AlgorithmBAcheck}) to $J$
using a reverse lexicographic order. If algorithm (\ref{AlgorithmBAcheck}) 
finishes without halting in step (\ref{Stephalt}), then Sullivant has 
proved that $J$ must generate $P$. If not, algorithm (\ref{AlgorithmBAcheck}) 
will halt with a binomial in $P\setminus J$.

Running the {\tt gbasis} command of {\tt GLATWALK} with respect
to the cost vector $-e_1$ and the grading $e_1 + \dots + e_{64}$ 
we compute a Gr\"obner basis of $J$ after converting the binomials
in the two files\footnote{
\tt http://math.berkeley.edu/\~{}seths/polyout.mac.gz} 
\footnote{\tt http://math.berkeley.edu/\~{}seths/polyout2.mac.gz}
containing $J$ into integer vector format. After computing a 
$15$-truncated Gr\"obner basis, {\tt gbasis} (in the incarnation of algorithm 
(\ref{AlgorithmBAcheck})) outputs the degree $14$ binomial
\begin{align*}
&x_{311} x_{221} x_{431} x_{212}^2 x_{122} x_{342} x_{113} x_{433} x_{243}
x_{424} x_{134} x_{334} x_{444}-\\
&x_{211} x_{421} x_{331} x_{112} x_{312}
x_{222} x_{242} x_{213} x_{133} x_{443} x_{124} x_{434}^2 x_{344}
\end{align*}
as a binomial in $\bar{J}\setminus J$ proving that $J$ does not
generate $P$ thereby answering Sullivant's computational commutative
algebra challenge. Running {\tt gbasis} in the above setting is not a
simple computation. In fact the $15$-truncated Gr\"obner basis 
of $J$ contains more than $300,000$ binomials and the whole
computation takes close to two days on most modern PCs.

Details and more information, including the relevant files for Sullivant's
challenge,  are located at {\tt http://home.imf.au.dk/niels/GLATWALK}.


\begin{thebibliography}{99}


\bibitem{CKR}

K. Caboara, M. Kreuzer, L. Robbiano. Effeciently computing minimal sets
of critical pairs. \emph{J.~Symbolic Computation}~{\bf 38} (2004), 1169--1190.

\bibitem{CLO}
Cox, Little and O'Shea. Ideals, Varieties and Algorithms. Undergraduate
Texts in Mathematics. Springer Verlag, 1992.

\bibitem{St}
B. Sturmfels. Gr\"obner Bases and Convex Polytopes, University Lecture Series
{\bf 8}, Amer. Math. Soc., Providence, RI, 1996.



\end{thebibliography}
\end{document}